\documentclass{amsart}
\usepackage[psamsfonts]{amssymb}
\usepackage{amsmath}
\usepackage{graphicx}
\usepackage[english]{babel}
\usepackage[all]{xy}
\usepackage{pifont}
\usepackage[pdftex]{hyperref}
\newtheorem{thm}{Teorema}[section]

\newtheorem{lem}[thm]{Lemma}

\theoremstyle{definition}

\numberwithin{equation}{section}

\begin{document}
\title{Energy of solenoidal vector fields on spherical domains}
\author{\textbf{Fabiano Brito \& Andr\'e Gomes}}%
\address{Instituto de Matem\'atica e Estat\'{\i}stica da Universidade de S\~ao Paulo}%
\email{gomes@ime.usp.br}%
\thanks{}%
\begin{abstract}
We present a ``boundary version" of a theorem about solenoidal unit vector fields with minimum energy on a spherical domain of an odd dimensional Euclidean sphere. 
\end{abstract}
\maketitle
\section{Introduction.}
Let $(M,g)$ be a closed, n-dimensional Riemannian manifold and $T^{1}M$ the unit tangent bundle of $M$ considered as a closed Riemannian manifold with the Sasaki
metric. Let $X:M\longrightarrow T^{1}M$ be a unit vector field defined on $M$, regarded as a smooth section of the unit tangent bundle $T^{1}M$. Using an orthonormal local frame $\left\{e_{1}, e_{2},\ldots,e_{n-1}, e_{n}=X\right\}$, the energy of the unit vector field $X$ is given by
\begin{eqnarray*}
	\mathcal{E}(X)=\frac{n}{2}\mathrm{vol}(M)+\frac{1}{2}\int_{M}\sum\limits_{a=1}^{n}\left\|\nabla_{e_{a}}X\right\|^{2}\nu_{_{M}}(g)
\end{eqnarray*}
The Hopf vector fields on $\mathbb{S}^{2k+1}$ are unit vector fields tangent to the classical Hopf
fibration $\mathbb{S}^{1}\hookrightarrow \mathbb{S}^{2k+1}$. The following theorem gives a characterization of Hopf flows as absolute minima of energy functional among all \textit{solenoidal} (that is, divergence free) unit vector fields on $\mathbb{S}^{2k+1}$

\normalfont 
\thm {[3] The Hopf vector fields has minimum energy among all solenoidal unit vector fields on the sphere $\mathbb{S}^{2k+1}$.}
\\
\\
\normalfont We prove in this paper the following boundary version for this Theorem:

\thm {Let $U$ be an open set of the $(2k+1)$-dimensional unit sphere $\mathbb{S}^{2k+1}$ and let $K\subset U$ be a connected $(2k+1)$-submanifold with boundary of the sphere $\mathbb{S}^{2k+1}$. Let $\vec{v}$ be a solenoidal unit vector field on $U$ which coincides with a Hopf flow $H$ along the boundary of K. Then 
\begin{eqnarray*}
	\mathcal{E}(\vec{v})\geq \left(\frac{2k+1}{2}+k\right)\mathrm{vol}(K)= \mathcal{E}(H)
\end{eqnarray*}
\normalfont
\\
\section{Preliminaries.}
Let $U\subset\mathbb{S}^{2k+1}$ be an open set of the unit sphere and let $K\subset U$ be a connected $(2k+1)$-submanifold with boundary of $\mathbb{S}^{2k+1}$. Let $H$ be a Hopf vector field on $\mathbb{S}^{2k+1}$ and let $\vec{v}$ be an unit vector field defined on $U$. We also consider the map $\varphi_{t}^{\vec{v}}:U\longrightarrow \mathbb{S}^{2k+1}(\sqrt{1+t^{2}})$ given by $\varphi_{t}^{\vec{v}}(x)=x+t\vec{v}(x)$. This map was introduced in [1] and [9].

\lem {For $t>0$ sufficiently small, the map $\varphi_{t}^{\vec{v}}$ is a diffeomorphism.}
\normalfont

\proof A simple application of the identity perturbation method $\square$
\\
\\
From now on, we assume that $t>0$ is small enough so that the map $\varphi_{t}^{\vec{v}}$ is a diffeomorphism. In order to find the Jacobian matrix of $\varphi_{t}^{\vec{v}}$, we define the unit vector field $\vec{u}$ on $\varphi_{t}^{\vec{v}}(U)\subset \mathbb{S}^{2k+1}(\sqrt{1+t^{2}})$ by
\begin{eqnarray*}
	\vec{u}(x):=\frac{1}{\sqrt{1+t^{2}}}\vec{v}(x)-\frac{t}{\sqrt{1+t^{2}}}x
\end{eqnarray*}
Using an adapted orthonormal frame $\left\{e_{1},\ldots, e_{2k},\vec{v}\right\}$ on a neighborhood $V$ of $U$, we obtain an adapted orthonormal frame on $\varphi_{t}^{\vec{v}}(V)$ given by $\left\{\bar{e}_{1},\ldots, \bar{e}_{2k},\vec{u}\right\}$, where $\bar{e}_{i}=e_{i}$ for all $i\in\left\{1,\ldots,2k\right\}$.
\\
\\
In this manner, we can write
\begin{eqnarray*}
	d\varphi_{t}^{\vec{v}}(e_{1})\!\!\!&=&\!\!\!\left\langle d\varphi_{t}^{\vec{v}}(e_{1}),e_{1}\right\rangle e_{1}+\ldots+\left\langle d\varphi_{t}^{\vec{v}}(e_{1}),e_{2k}\right\rangle e_{2k}+\left\langle d\varphi_{t}^{\vec{v}}(e_{1}),\vec{u}\right\rangle \vec{u}\\
	d\varphi_{t}^{\vec{v}}(e_{2})\!\!\!&=&\!\!\!\left\langle d\varphi_{t}^{\vec{v}}(e_{2}),e_{1}\right\rangle e_{1}+\ldots+\left\langle d\varphi_{t}^{\vec{v}}(e_{2}),e_{2k}\right\rangle e_{2k}+\left\langle d\varphi_{t}^{\vec{v}}(e_{2}),\vec{u}\right\rangle \vec{u}\\
	\vdots
	\\
	d\varphi_{t}^{\vec{v}}(e_{2k})\!\!\!&=&\!\!\!\left\langle d\varphi_{t}^{\vec{v}}(e_{2k}),e_{1}\right\rangle e_{1}+\ldots+\left\langle d\varphi_{t}^{\vec{v}}(e_{2k}),e_{2k}\right\rangle e_{2k}+\left\langle d\varphi_{t}^{\vec{v}}(e_{2k}),\vec{u}\right\rangle \vec{u}\\
	d\varphi_{t}^{\vec{v}}(\vec{v})\!\!\!&=&\!\!\!\left\langle d\varphi_{t}^{\vec{v}}(\vec{v}),e_{1}\right\rangle e_{1}+\ldots+\left\langle d\varphi_{t}^{\vec{v}}(\vec{v}),e_{2k}\right\rangle e_{2k}+\left\langle d\varphi_{t}^{\vec{v}}(\vec{v}),\vec{u}\right\rangle \vec{u}
\end{eqnarray*}
Now, by Gauss' equation of the trivial immersion $\mathbb{S}^{2k+1}\hookrightarrow \mathbb{R}^{2k+2}$, we have 
\begin{eqnarray*}
	\tilde{\nabla}_{Y}\vec{v}=d\vec{v}(Y)=\nabla_{Y}\vec{v}-\left\langle \vec{v},Y\right\rangle x
\end{eqnarray*}
for every vector field $Y$ on $\mathbb{S}^{2k+1}$, and then
\begin{eqnarray*}
	\left\langle d\varphi_{t}^{\vec{v}}(e_{1}),e_{1}\right\rangle = \left\langle e_{1}+td\vec{v}(e_{1}),e_{1}\right\rangle = 1 + t\left\langle \nabla_{e_{1}}\vec{v},e_{1}\right\rangle
\end{eqnarray*}
Analogously, we can conclude that 
\begin{eqnarray*}
\left\langle d\varphi_{t}^{\vec{v}}(e_{i}),e_{i}\right\rangle \!\!\!&=&\!\!\!1+t\left\langle \nabla_{e_{i}}\vec{v},e_{i}\right\rangle, \ \forall i\in\left\{1,\ldots,2k\right\}\\
\left\langle d\varphi_{t}^{\vec{v}}(e_{i}),e_{j}\right\rangle \!\!\!&=&\!\!\! t\left\langle \nabla_{e_{i}}\vec{v},e_{j}\right\rangle, \ \forall i,j \in \left\{1,\ldots,2k\right\},\ (i\neq j)\\
\left\langle d\varphi_{t}^{\vec{v}}(e_{i}),\vec{u}\right\rangle \!\!\!&=&\!\!\! 0, \ \forall i\in \left\{1,\ldots,2k\right\}\\
\left\langle d\varphi_{t}^{\vec{v}}(\vec{v}),\vec{u}\right\rangle \!\!\!&=&\!\!\! \sqrt{1+t^{2}}
\end{eqnarray*}
By employing the notation $h_{ij}(\vec{v}):=\left\langle \nabla_{e_{i}}\vec{v},e_{j}\right\rangle$ (where $i,j\in\left\{1,\ldots,2k\right\}$), we can express the determinant of the Jacobian matrix of $\varphi_{t}^{\vec{v}}$ in the form 
\begin{eqnarray*}
	\det(d\varphi_{t}^{\vec{v}})=\sqrt{1+t^{2}}(1+\sum\limits_{i=1}^{2k}\sigma_{i}(\vec{v})t^{2})
\end{eqnarray*}
where, by definition, the functions $\sigma_{i}$ are the $i$-symmetric functions of the $h_{ij}$. For instance, if $k=1$, we have
\begin{eqnarray*}
\sigma_{1}(\vec{v})\!\!\!&:=&\!\!\!h_{11}(\vec{v})+h_{22}(\vec{v})\\ \sigma_{2}(\vec{v})\!\!\!&:=&\!\!\!h_{11}(\vec{v})h_{22}(\vec{v})-h_{12}(\vec{v})h_{21}(\vec{v})
\end{eqnarray*}

\section{Proof of the Theorem.}
The energy of the vector field $\vec{v}$ (on $K$) is given by 
\begin{eqnarray*}	
\mathcal{E}(\vec{v}):=\frac{1}{2}\int_{K}\left\|d\vec{v}\right\|^{2}=\frac{2k+1}{2}\mathrm{vol}(K)+\frac{1}{2}\int_{K}\left\|\nabla \vec{v}\right\|^{2}
\end{eqnarray*}
Using the notation above, we have
\begin{eqnarray*}	\mathcal{E}(\vec{v})=\frac{2k+1}{2}\mathrm{vol}(K)+\frac{1}{2}\int_{K}[\sum\limits_{i,j=1}^{2k}(h_{ij}(\vec{v}))^{2}+\sum\limits_{i=1}^{2k}(\left\langle\nabla_{\vec{v}}\vec{v},e_{i}\right\rangle)^{2}]
\end{eqnarray*}
and then
\begin{eqnarray}
\mathcal{E}(\vec{v})\geq \frac{2k+1}{2}\mathrm{vol}(K)+\frac{1}{2}\int_{K}\sum\limits_{i,j=1}^{2k}(h_{ij}(\vec{v}))^{2}
\end{eqnarray}
Now observe that
\begin{eqnarray}
\sum\limits_{i<j}(h_{ii}-h_{jj})^{2}=(2k-1)\sum\limits_{i}h_{ii}^{2}-2\sum\limits_{i<j}h_{ii}h_{jj}
\end{eqnarray}
and as $\vec{v}$ is a solenoidal vector field
\begin{eqnarray}
	0=[div(\vec{v})]^{2}=[\sigma_{1}(\vec{v})]^{2}=(\sum\limits_{i}h_{ii})^{2}=\sum\limits_{i}h_{ii}^{2}+2\sum\limits_{i<j}h_{ii}h_{jj}
\end{eqnarray}
in other words
\begin{eqnarray}
	-2\sum\limits_{i<j}h_{ii}h_{jj}=\sum\limits_{i}h_{ii}^{2}
\end{eqnarray}
Substituting equation (3.4) in (3.2) we obtain
\begin{eqnarray}
\sum\limits_{i<j}(h_{ii}-h_{jj})^{2}=-4k\sum\limits_{i<j}h_{ii}h_{jj}
\end{eqnarray}
Further, we also have the following equation
\begin{eqnarray}
\sum\limits_{i<j}(h_{ij}+h_{ji})^{2}=\sum\limits_{i\neq j}h_{ij}^{2}+2\sum\limits_{i<j}h_{ij}h_{ji}
\end{eqnarray}
and then
\begin{eqnarray}
2k\sum\limits_{i<j}(h_{ij}+h_{ji})^{2}=2k\sum\limits_{i\neq j}h_{ij}^{2}+4k\sum\limits_{i<j}h_{ij}h_{ji}
\end{eqnarray}
Adding equations (3.5) and (3.7), we have
\begin{eqnarray}
\sum\limits_{i\neq j}h_{ij}^{2}\geq 2\sigma_{2}
\end{eqnarray}
and
\begin{eqnarray}
\sum\limits_{i,j=1}^{2k}h_{ij}^{2}=\sum\limits_{i}h_{ii}^{2}+\sum\limits_{i\neq j}h_{ij}^{2}\geq 2\sigma_{2}
\end{eqnarray}
Using the inequalities (3.1) and (3.9), we find 
\begin{eqnarray}
\mathcal{E}(\vec{v})&\geq&\!\!\! \frac{2k+1}{2}\mathrm{vol}(K)+\int_{K}\sigma_{2}(\vec{v})
\end{eqnarray}
\\
On the other hand, by change of variables theorem, we obtain
\begin{eqnarray*}	\mathrm{vol}[\varphi_{t}^{H}(K)]=\int_{K}\sqrt{1+t^{2}}(1+\sum\limits_{i=1}^{2k}\sigma_{i}(H)t^{i})
\end{eqnarray*}
By a straightforward computation shown in [4], we have $\sigma_{i}(H)=\eta_{i}$ for all index $i\in\left\{1,\ldots,2k\right\}$, where the numbers $\eta_{i}$ are defined by
\[
  \eta_{i} = \left\{ 
  \begin{array}{l l}
    \!\!{k\choose i/2} & \quad i\!f \  i \  is \  even\\
    \\
    0 & \quad i\!f \  i \  is \   odd \\
  \end{array} \right.
\]
\\
We know that the vector fields $\vec{v}$ and $H$ are the same on $\partial K$. Thus, $\varphi_{t}^{\vec{v}}(K)$ and $\varphi_{t}^{H}(K)$ are $(2k+1)$-submanifolds of $\mathbb{S}^{2k+1}(\sqrt{1+t^{2}})$ with the same boundary. We claim that $\varphi_{t}^{\vec{v}}(K)=\varphi_{t}^{H}(K)$ for all $t$ sufficiently small. In fact, if $p$ is an interior point of $K$, $\lim\limits_{t\rightarrow 0}\varphi_{t}^{\vec{v}}(p)=\lim\limits_{t\rightarrow 0}\varphi_{t}^{H}(p)=p$ and then we have necessarily 
\begin{eqnarray*}
	\mathrm{vol}[\varphi_{t}^{\vec{v}}(K)]=\mathrm{vol}[\varphi_{t}^{H}(K)]
\end{eqnarray*}
for all $t$ sufficiently small, or equivalently,
\begin{eqnarray*}	
\int_{K}\sqrt{1+t^{2}}(1+\sum\limits_{i=1}^{2k}\sigma_{i}(\vec{v})t^{i})=\int_{K}\sqrt{1+t^{2}}(1+\sum\limits_{i=1}^{2k}\eta_{i}t^{i})
\end{eqnarray*}
for all $t>0$ sufficiently small. Consequently, after cancelling the factor $\sqrt{1+t^{2}}$ and rearranging the terms, we obtain 
\begin{eqnarray*}
\left(\int_{K}[\sigma_{1}(\vec{v})-\eta_{1}]\right)t+\left(\int_{K}[\sigma_{2}(\vec{v})-\eta_{2}]\right)t^{2}+\ldots+\left(\int_{K}[\sigma_{2k}(\vec{v})-\eta_{2k}]\right)t^{2k}= 0
\end{eqnarray*}
for all sufficiently small $t$. By identity of polynomials, we conclude
\begin{eqnarray}
	\int_{K}\sigma_{i}(\vec{v})=\int_{K}\eta_{i}=\eta_{i}\mathrm{vol}(K),\ \ \forall i\in\left\{1,\ldots,2k\right\}
\end{eqnarray}
and then, using the inequality (3.10) and the equality (3.11) (for $i=2$), we have
\begin{eqnarray*}
	\mathcal{E}(\vec{v})\geq \frac{2k+1}{2}\mathrm{vol}(K)+\eta_{2}\mathrm{vol}(K)=\left(\frac{2k+1}{2}+k\right)\mathrm{vol}(K)
\end{eqnarray*}


\begin{thebibliography}{99}

\bibitem[1]{01}D. Asimov, \textit{Average Gaussian curvature of leaves of foliations}, Bull. Amer. Math. Soc. 84, (1978), 131-133
\bibitem[2]{02}F. G. B. Brito, \textit{Total bending of flows with mean curvature correction}, Diff. Geom. Appl. 12, (2000), 157-163.

\bibitem[3]{03}F. G. Brito, M. Salvai, \textit{Solenoidal unit vector fields with minimum energy}, Osaka J. Math. 41, (2004), 533-544.

\bibitem[4]{04}P. M. Chac\'on, A. M. Naveira, F. B. Brito, \textit{On the volume of unit vector fields on spaces of constant sectional curvature}, Comment. Math. Helv. 79 (2004), no. 2, 300-316. 
\bibitem[5]{05}P. M. Chac\'on, G. S. Nunes, \textit{Energy and topology of singular unit vector fields on $\mathbb{S}^{3}$}, Pacific J. Math. 231 (2007), no. 1, 27-34. 
\bibitem[6]{06}P. M. Chac\'on, A. M. Naveira, J. M. Weston, \textit{On the energy of distributions, with application to the quaternionic Hopf fibrations}, Monatsh. Math. 133 (2001), no. 4, 281-294. 
\bibitem[7]{07}O. Gil-Medrano, A. Hurtado, \textit{Volume, energy and generalized energy of unit vector fields on Berger spheres: stability of Hopf vector fields}, Proc. Roy. Soc. Edinburgh Sect. A 135 (2005), no. 4, 789-813.
\bibitem[8]{08}O. Gil-Medrano, E. Llinares-Fuster, \textit{Minimal unit vector fields}, Tohoku Math. J. (2) 54 (2002), no. 1, 71-84.
\bibitem[9]{09}J. Milnor, \textit{Analytic proofs of the ``hairy ball theorem" and the Brouwer fixed-point theorem}, Amer. Math. Monthly 85, (1978), no. 7, 521-524.
\bibitem[10]{10}L. Vanhecke, E. Boeckx, \textit{Harmonic and minimal radial vector fields}, Acta Math. Hungar. 90 (2001), no. 4, 317-331. 
\bibitem[11]{11}L. Vanhecke, E. Boeckx, J. C. Gonz\'alez-D\'avila, \textit{Energy of radial vector fields on compact rank one symmetric spaces}, Ann. Global Anal. Geom. 23 (2003), no. 1, 29-52. 

\end{thebibliography}
\end{document}